\def\draftdate{August 22, 2007}

\documentclass{amsart}
\usepackage{xy}
\xyoption{arrow}
\xyoption{curve}
\xyoption{matrix}
\xyoption{cmtip}
\SelectTips{cm}{}

\newcommand{\ssdot}{\bullet}
\newcommand{\subdot}{_\ssdot}

\newcommand{\Sdot}[1][\ssdot]{S_{#1}}
\newcommand{\Spdot}[1][\ssdot]{S'_{#1}}
\newcommand{\phat}{^{\scriptscriptstyle\wedge}_{p}}
\newcommand{\ang}[1][{n}]{\langle #1\rangle}
\newcommand{\cU}[1][{n+1}]{\aU_{#1}}

\let\sma\wedge
\renewcommand{\to}{\mathchoice{\longrightarrow}{\rightarrow}{\rightarrow}{\rightarrow}}

\newcommand{\bZ}{{\mathbb{Z}}}
\let\catsymbfont\mathcal 

\newcommand{\aB}{{\catsymbfont{B}}}
\newcommand{\aC}{{\catsymbfont{C}}}
\newcommand{\aD}{{\catsymbfont{D}}}
\newcommand{\aE}{{\catsymbfont{E}}}

\newcommand{\aM}{{\catsymbfont{M}}}
\newcommand{\aP}{{\catsymbfont{P}}}
\newcommand{\aU}{{\catsymbfont{U}}}

\def\quickop#1{\expandafter\DeclareMathOperator\csname
#1\endcsname{#1}}
\quickop{id}\quickop{Id}\quickop{colim}\quickop{hocolim}\quickop{op}
\quickop{co}\quickop{Ar}\quickop{sing}\quickop{Hom}

\newtheorem{thm}{Theorem}[section]

\newtheorem{lem}[thm]{Lemma}
\newtheorem{prop}[thm]{Proposition}
\newtheorem*{arc}{Rognes Conjecture}
\newtheorem*{main-l}{Localization Theorem}
\newtheorem*{main-d}{D\'{e}vissage Theorem}
\theoremstyle{definition}
\newtheorem{defn}[thm]{Definition}
\theoremstyle{remark}
\newtheorem{rem}[thm]{Remark}

\numberwithin{equation}{section}
\makeatletter\let\c@equation\c@thm\makeatother

\newcommand{\term}[1]{\textit{#1}}

\newdir{ >}{{}*!/-5pt/\dir{>}}

\begin{document}

\title[The algebraic $K$-theory of topological $K$-theory]%
{The localization sequence for the algebraic $K$-theory of
topological $K$-theory}

\author{Andrew J. Blumberg}
\address{Department of Mathematics, Stanford  University,
Stanford, CA 94305}
\email{blumberg@math.stanford.edu}
\thanks{The first author was supported in part by an NSF postdoctoral
fellowship}

\author{Michael A. Mandell}
\address{Department of Mathematics, Indiana University,
Bloomington, IN \ 47405}
\email{mmandell@indiana.edu}
\thanks{The second author was supported in part by NSF grant DMS-0504069}

\date{\draftdate}

\begin{abstract}
We verify a conjecture of Rognes by establishing a localization
cofiber sequence of spectra
$K(\mathbb{Z})\to K(ku)\to K(KU) \to \Sigma K(\mathbb{Z})$
for the algebraic $K$-theory of topological $K$-theory.  We deduce
the existence of this sequence as a consequence of a d\'{e}vissage theorem
identifying the $K$-theory of the Waldhausen category of
finitely generated finite stage Postnikov towers of modules over a
connective $A_\infty$ ring spectrum $R$ with the Quillen $K$-theory of
the abelian category of finitely generated $\pi_{0}R$-modules.
\end{abstract}

\maketitle

\section*{Introduction}

The algebraic $K$-theory of ring spectra is intimately related to
the geometry of high-dimensional manifolds.  In a series of papers in
the 1970's and 1980's, Waldhausen established the deep connection
between the $K$-theory of the sphere spectrum and stable
pseudo-isotopy theory.  The sphere spectrum has an infinite filtration
called the chromatic filtration that forms a tower at each prime
$p$.  The layers in the chromatic tower capture periodic phenomena in
stable homotopy theory, corresponding to the Morava $K$-theory
``fields''.  The chromatic viewpoint has organized the understanding of stable
homotopy theory over the last thirty years.

Applying algebraic $K$-theory to the chromatic tower of the sphere
spectrum leads to an analogous ``chromatic tower'' of algebraic
$K$-theory spectra.  Waldhausen conjectured that this tower converges
in the homotopy inverse limit to the algebraic $K$-theory of the
$p$-local sphere spectrum \cite{Wchromatic}.  McClure and Staffeldt
verified a closely related connective variant of this conjecture
\cite{McClureStaffeldt}.  The bottom layer of the chromatic tower is  
$K(\mathbb{Q})$, the algebraic $K$-theory of the rational numbers,
which is intimately connected to questions in arithmetic.  Thus, the
$K$-theory chromatic tower can be viewed as an interpolation from
arithmetic to geometry.

In \cite{AusoniRognes}, Ausoni and Rognes describe an ambitious
program for analyzing the layers in the $p$-complete version of
Waldhausen's $K$-theory chromatic tower.  The analysis is based on
descent conjectures coming from Rognes' Galois theory of $S$-algebras
\cite{RognesGalois}, which relate the layers of the tower to the
$K$-theory of the Morava $E$-theory ring spectra $E_{n}$ and to the
$K$-theory of the $p$-completed Johnson--Wilson ring spectra
$E(n)\phat$.  The spectrum $E(n)$ is not connective, and no tools
exist for computing the $K$-theory of non-connective ring spectra.
However, $E(n)$ is formed from the connective spectrum $BP\ang$ by
inverting the element $v_{n}$ in $\pi_{*}(BP\ang)$, and Rognes
conjectured the following localization cofiber sequence:
\begin{arc}
The transfer
map $K(BP\ang[n-1]\phat)\to K(BP\ang\phat)$ and the canonical map
$K(BP\ang\phat)\to K(E(n)\phat)$ fit into a cofiber
sequence in the stable category
\[ K(BP\ang[n-1]\phat)\to K(BP\ang\phat)\to K(E(n)\phat)
\to  \Sigma K(BP\ang[n-1]\phat).
\]
\end{arc}
In the case $n=0$, the statement is an old theorem of Quillen
\cite{QuillenAK}, the localization
sequence $K(\mathbb{Z}/p)\to K(\mathbb{Z}\phat)\to
K(\mathbb{Q}\phat)$.  The conjecture implies a long exact
sequence of homotopy groups: 
\begin{multline*} \dotsb \to
K_{q}(BP\ang[n-1]\phat)\to K_{q}(BP\ang\phat)\to K_{q}(E(n)\phat) \to
\dotsb \\
\dotsb \to
K_{0}(BP\ang[n-1]\phat)\to K_{0}(BP\ang\phat)\to K_{0}(E(n)\phat)\to 0 \end{multline*}
This in principle allows the
computation of the algebraic $K$-groups of $E(n)\phat$ from those of
$BP\ang\phat$ and $BP\ang[n-1]\phat$.  Because the $BP\ang$ are
connective spectra, their $K$-theory can be studied using $TC$, and
the main purpose of \cite{AusoniRognes} is to provide tools suitable for
computing $TC(BP\ang\phat)$ and in particular to obtain an explicit
computation of the $V(1)$-homotopy in the case $n=1$.

The case $n=1$ is the first non-classical case of the conjecture and
is of particular interest.  The spectra $E(1)$ and $BP\ang[1]$ are
often denoted $L$ and $\ell$ and are closely related to topological
$K$-theory.  Here $L = E(1)$ is the $p$-local Adams summand of
(complex periodic) topological $K$-theory, $KU$, and $\ell =
BP\ang[1]$ is the $p$-local Adams summand of connective topological
$K$-theory, $ku$.  The Rognes conjecture then admits a ``global''
version in this case, namely that the transfer map $K(\mathbb{Z})\to
K(ku)$ and canonical map $K(ku)\to K(KU)$ fit into a cofiber sequence
in the stable category 
\[
K(\mathbb{Z})\to K(ku)\to K(KU) \to \Sigma K(\mathbb{Z}).
\]
This version of the
conjecture appears in \cite{BaasDundasRognes}, where the relationship
between the algebraic $K$-theory of topological $K$-theory, elliptic
cohomology, and the category of $2$-vector bundles is discussed.

The purpose of this paper is to prove both the local and global
versions of the Rognes conjecture in the case $n=1$.
Specifically, we prove the following theorem:

\begin{main-l}
There are connecting maps, which together with the transfer maps and the
canonical maps, make the sequences
\begin{gather*}
K(\mathbb{Z}\phat)\to K(\ell\phat)\to K(L\phat)\to \Sigma K(\mathbb{Z}\phat)\\
K(\mathbb{Z}_{(p)})\to K(\ell)\to K(L)\to \Sigma K(\mathbb{Z}_{(p)})\\
K(\mathbb{Z})\to K(ku)\to K(KU) \to \Sigma K(\mathbb{Z})
\end{gather*}
cofiber sequences in the stable category.
\end{main-l}

Hesselholt observed that the Ausoni-Rognes and Ausoni calculations of $THH(\ell)$
and $THH(ku)$ in \cite{AusoniRognes, Ausoni} are explained by the
existence of a $THH$ version of this localization sequence along
with a conjecture about the behavior of $THH$ for ``tamely ramified''
extensions of ring spectra.  A precise formulation requires a
construction of $THH$ for Waldhausen categories.  We will explore
this more fully in a forthcoming paper.

The Localization Theorem above is actually a consequence of a
``d\'{e}vissage'' theorem for finitely generated finite stage Postnikov
towers.  For an $S$-algebra ($A_{\infty}$ ring spectrum) $R$, let
$\aP_{R}$ denote the full subcategory of left $R$-modules that have
only finitely many non-zero homotopy groups, all of which are finitely
generated over $\pi_{0}R$.  When $R$ is connective and $\pi_{0}R$ is
left Noetherian, this category has an associated Waldhausen $K$-theory
spectrum.  Restricting to the subcategory of $S$-algebras with
morphisms the maps $R\to R'$ for which $\pi_{0}R'$ is
finitely generated as a left $\pi_{0}R$-module, we regard $K(\aP_{(-)})$ as
a contravariant functor $K'$ to the stable category.

We use the notation $K'$ because of the close connection with
Quillen's $K'$-theory, the $K$-theory of the exact category of finitely generated
left modules over a left Noetherian ring.  The analogous $K$-theory
spectrum for chain complexes over the ring $\pi_{0}R$ is well-known to
be equivalent to $K'(\pi_{0}R)$ \cite[1.11.7]{TTGrothFest}.  The
following theorem is the main result of this paper.

\begin{main-d}
Let $R$ be a connective $S$-algebra ($A_{\infty}$ ring spectrum) with
$\pi_{0}R$ left Noetherian.  Then there is a natural isomorphism in
the stable category $K'(\pi_{0}R)\to K'(R)$, where $K'(\pi_{0}R)$ is
Quillen's $K$-theory of the exact category of finitely generated left
$\pi_{0}R$-modules, and $K'(R)$ is the Waldhausen $K$-theory of the
category of finitely generated finite stage Postnikov towers of left
$R$-modules, $\aP_{R}$.
\end{main-d}

A longstanding open problem first posed explicitly by Thomason
and Trobaugh \cite[1.11.1]{TTGrothFest} is to develop a general d\'{e}vissage
theorem for Waldhausen categories that specializes to Quillen's
d\'{e}vissage theorem when applied to the category of bounded chain
complexes on an abelian category.  We regard the theorem described in
this paper as a step towards a solution to this problem.

\bigskip

The authors would like to thank the Institut Mittag--Leffler for
hospitality while writing this paper.

\section{The Main Argument}\label{secmainarg}

In this section, we outline the proof of the D\'{e}vissage Theorem in
terms of a number of easily stated results proved in later sections;
we then deduce the Localization Theorem from the D\'{e}vissage Theorem.
Although we assume some familiarity with the basics of Waldhausen
$K$-theory, we review the standard definitions and constructions
of Waldhausen \cite{Wald} as needed.  We begin with some
technical conventions and a precise description of the Waldhausen
categories we use.

Throughout this paper, $R$ denotes a connective $S$-algebra with
$\pi_{0}R$ left Noetherian.  We work in the context of EKMM $S$-modules,
$S$-algebras, and $R$-modules \cite{EKMM}.  Since other contexts for
the foundations of a modern category of spectra lead to equivalent
$K$-theory spectra, presumably the arguments presented here could be
adjusted to these contexts, but the EKMM categories have certain
technical advantages that we exploit (see for example the proof of
Lemma~\ref{lemapplymainlem}) and that affect the precise form of the
statements below.  

The input for Waldhausen $K$-theory is a (small) category together
with a subcategory of ``weak equivalences'' and a subcategory of
``w-cofibrations'' satisfying certain properties.  Let $\aP_{R}$
denote the full subcategory of left $R$-modules that have only
finitely many non-zero homotopy groups, all of which are finitely
generated over $\pi_{0}R$.  Although $\aP_{R}$ is not a small
category, we can still construct from it a $K$-theory spectrum that is
``homotopically small'', and we can find a small category with
equivalent $K$-theory by restricting the sets allowed in the
underlying spaces of the underlying prespectra; see
Remark~\ref{remsmall} below for details.  In what follows, let $\aP$
denote $\aP_{R}$, or at the reader's preference, the small category
$\aP_{\kappa}$ (for $\kappa$ large) discussed in
Remark~\ref{remsmall}.

We make $\aP$ a Waldhausen category by taking the weak equivalences to
be the usual weak equivalences (the maps that induce isomorphisms on all
homotopy groups) and the w-cofibrations to be the Hurewicz
cofibrations (the maps satisfying the homotopy extension property in the
category of left $R$-modules).  Specifically, a map $i \colon A \to
X$ is a Hurewicz cofibration if and only if the inclusion of the mapping
cylinder $Mi = X \cup_i (A \sma I_+)$ in the cylinder $X\sma I_{+}$
has a retraction in the category of $R$-modules.  Some easy consequences
of this definition are:
\begin{itemize}
\item The initial map $*\to X$ is a Hurewicz cofibration for any $R$-module $X$,
\item Hurewicz cofibrations are preserved by cobase change (pushout), and
\item For a map of $S$-algebras $R\to R'$, the forgetful functor from
$R'$-modules to $R$-modules preserves Hurewicz cofibrations.
\end{itemize}
Use of this type of cofibration was a key tool in \cite{EKMM} for
keeping homotopical control; a key fact is that pushouts along
Hurewicz cofibrations preserve weak equivalences \cite[I.6.5]{EKMM},
and in particular Waldhausen's Gluing Lemma \cite[1.2]{Wald} holds.

Since we are thinking of $\aP$ in terms of finite Postnikov stages, it
makes more sense philosophically to work with fibrations: $\aP^{\op}$
forms a Waldhausen category with weak equivalences the maps opposite
to the usual weak equivalences and with w-cofibrations the maps
opposite to the Hurewicz fibrations (maps satisfying the covering
homotopy property).  Although not strictly necessary for the
D\'{e}vissage Theorem, the following theorem proved in
Section~\ref{seccoffib} straightens out this discrepancy.

\begin{thm}\label{thmcoffib}
The spectra $K(\aP)$ and $K(\aP^{\op})$ are weakly equivalent.
\end{thm}

This result is essentially a consequence of the fact that $\aP$ is a
stable category, and in particular follows from the observation that
homotopy cocartesian and homotopy cartesian squares coincide in $\aP$.
As part of the technical machinery employed in the proof of the
D\'{e}vissage Theorem, we introduce a variant of Waldhausen's $\Sdot$
construction, using homotopy cocartesian squares where 
pushout squares along w-cofibrations are used in the $\Sdot$
construction.  We denote this construction as $\Spdot$.
Theorem~\ref{thmcoffib} follows from the 
comparison of the $\Spdot$ construction with the $\Sdot$ construction
in Section~\ref{seccoffib}.

Let $\aP_{m}^{n}$ for $m\leq n$ denote the full subcategory of $\aP$
consisting of those $R$-modules whose homotopy groups $\pi_{q}$ are
zero for $q>n$ or $q<m$.  In this notation, we permit $m=-\infty$ and/or
$n=\infty$, so $\aP=\aP_{-\infty}^{\infty}$.  Define a 
w-cofibration in $\aP_{m}^{n}$ to be a Hurewicz cofibration whose
cofiber is still in $\aP_{m}^{n}$, or equivalently, a Hurewicz
cofibration inducing an injection on $\pi_{n}$.  This definition makes
$\aP_{m}^{n}$ into a Waldhausen category with the usual weak
equivalences; it is a ``Waldhausen subcategory'' of $\aP$
\cite[1.2]{Wald}.  We prove the following theorem in
Sections~\ref{secbigthm} and~\ref{secpizero}. 

\begin{thm}\label{thmbigthm}
The inclusion $\aP_{0}^{0}\to \aP$ induces a weak equivalence of
$K$-theory spectra.
\end{thm}

Let $\aE$ denote the exact category of finitely generated left
$\pi_{0}R$-modules (or a skeleton, to obtain a small category).  This
becomes a Waldhausen category with weak equivalences the isomorphisms
and w-cofibrations the injections.  Waldhausen's ``$\Sdot=Q$'' Theorem
\cite[1.9]{Wald} identifies the Waldhausen $K$-theory $K(\aE)$ as
$K'(\pi_{0}R)$.  The functor $\pi_{0}\colon \aP_{0}^{0}\to \aE$ is an
``exact functor'' of Waldhausen categories: It preserves weak
equivalences, w-cofibrations, and pushouts along w-cofibrations.  It
follows that $\pi_{0}$ induces a map from $K(\aP_{0}^{0})$ to
$K(\aE)\simeq K'(\pi_{0}R)$.  We prove the following theorem in
Section~\ref{secpizero}.

\begin{thm}\label{thmpizero}
The functor $\pi_{0}$ induces a weak equivalence $K(\aP^{0}_{0})\to K(\aE)$.
\end{thm}

Theorem~\ref{thmpizero} and Theorem~\ref{thmbigthm} together imply the
D\'{e}vissage Theorem.  For $K$-theoretic reasons, we should regard
$K'(\pi_{0}R)\to K(\aP)$ to be the natural direction of the composite
zigzag, as this is compatible with the forgetful functor
$\aP_{H\pi_{0}R}\to \aP_{R}$ (induced by pullback along the map $R \to
H\pi_{0}R$).  Although it appears feasible to construct directly a map
of spectra $K'(\pi_{0}R)\to K(\aP)$ using a version of the
Eilenberg--Mac\,Lane bar construction, the technical work required
would be unrelated to the arguments in the rest of this paper, and so
we have not pursued it.

Next we deduce the Localization Theorem from the D\'{e}vissage Theorem.
Let $R$ be one of $ku$, $\ell$, or $\ell\phat$, and let $\beta$ denote
the appropriate Bott element in $\pi_{*}R$ in degree $2$ or $2p-2$.
Then $R[\beta^{-1}]$ is $KU$, $L$, or $L\phat$ respectively.  For
convenience, let $Z$ denote $\pi_{0}R$; so $Z=\bZ$, $\bZ_{(p)}$, or
$\bZ\phat$ in the respective cases.  Then for $A=HZ$, $R$, or
$R[\beta^{-1}]$, let $\aC_{A}$ be the category of finite cell
$A$-modules, which we regard as a Waldhausen category with
w-cofibrations the maps that are isomorphic to inclusions of cell
subcomplexes and with weak equivalences the usual weak equivalences,
that is, the maps that induce isomorphisms on homotopy groups.  As
explained in \cite{EKMM}, the Waldhausen $K$-theory of $\aC_{A}$ is
the algebraic $K$-theory of $A$.  When $A=HZ$, the $K$-theory spectrum
$K(\aC_{A})$ is equivalent to Quillen's $K$-theory of the ring $Z$
\cite[VI.4.3]{EKMM}.  In the case of $A=HZ$ or $A=R$, 
the $K$-theory spectrum $K(\aC_{A})$ is equivalent to two other
reasonable versions of the $K$-theory of $A$: the $K$-theory
spectrum defined via the plus construction of ``$BGL(A)$''
\cite[VI\S7]{EKMM} and the $K$-theory spectrum defined via the
permutative category of wedges of sphere $A$-modules
\cite[VI\S6]{EKMM}.  However, when $A=R[\beta^{-1}]$ the construction
$K(\aC_{A})$ and its variants are essentially the only known ways to
define $K(A)$.

Let $\aC_{R}[\beta^{-1}]$ denote the Waldhausen category whose
underlying category is $\aC_{R}$ and whose w-cofibrations are the maps
that are isomorphic to inclusions of cell subcomplexes, but whose weak
equivalences are the maps that induce isomorphisms on homotopy groups
after inverting $\beta$.  The identity functor $\aC_{R}\to
\aC_{R}[\beta^{-1}]$ is then an exact functor because the weak
equivalences in $\aC_{R}$ are in particular weak equivalences in
$\aC_{R}[\beta^{-1}]$.  Let $\aC_{R}^{\beta}$ be the full subcategory
of $\aC_{R}$ of objects that are weakly equivalent in
$\aC_{R}[\beta^{-1}]$ to the trivial object $*$, that is, the finite
cell $R$-modules whose homotopy groups become zero after inverting
$\beta$.  In other words, we consider $\aC_{R}$ with two
subcategories of weak equivalences, one of which is
coarser than the other, and the category of acyclic objects for the
coarse weak equivalences.  This is the situation in which Waldhausen's
``Fibration Theorem'' \cite[1.6.4]{Wald} applies.  The conclusion in
this case, restated in terms of the stable category, is the following
proposition.

\begin{prop}\label{propwaldfib}
There is a connecting map (of spectra) $K(\aC_{R}[\beta^{-1}])\to
\Sigma K(\aC_{R}^{\beta})$, which together with the canonical maps
makes the sequence
\[
K(\aC_{R}^{\beta})\to K(R)\to K(\aC_{R}[\beta^{-1}])\to \Sigma K(\aC_{R}^{\beta})
\]
a cofiber sequence in the stable category.
\end{prop}

For the proof of the Localization Theorem, we need to identify 
$K(\aC_{R}^{\beta})$ as $K(\pi_{0}R)$, $K(\aC_{R}[\beta^{-1}])$ as
$K(R[\beta^{-1}])$ and the maps as the transfer map and the canonical
map.   

We begin with the comparison of $K(\aC_{R}[\beta^{-1}])$ with
$K(R[\beta^{-1}])$.  The localization functor $R[\beta^{-1}] \sma_{R}
(-)$ restricts to an exact functor of Waldhausen categories
$\aC_{R}[\beta^{-1}] \to \aC_{R[\beta^{-1}]}$.  Waldhausen developed
a general tool for showing that an exact functor induces a weak
equivalence of $K$-theory spectra, called the ``Approximation
Theorem'' \cite[1.6.7]{Wald}.  In the case of a telescopic
localization functor like this one, the argument of \cite{Wchromatic}
adapts to show that the hypotheses of Waldhausen's Approximation
Theorem are satisfied.  This then proves the following proposition.

\begin{prop}\label{propidentright}
The canonical map $K(R)\to K(R[\beta^{-1}])$ factors as the map
$K(R)\to K(\aC_{R}[\beta^{-1}])$ in Proposition~\ref{propwaldfib} and
a weak equivalence $K(\aC_{R}[\beta^{-1}])\to K(R[\beta^{-1}])$.
\end{prop}

Next we move on to the comparison of $K(\aC_{R}^{\beta})$ with
$K(\pi_{0}R)$.  Since finite cell $R$-modules that are $\beta$-torsion
have finitely generated homotopy groups concentrated in a finite
range, $\aC_{R}^{\beta}$ is a subcategory of $\aP_{R}$.  Cellular
inclusions are Hurewicz cofibrations, so the inclusion
$\aC_{R}^{\beta}\to \aP_{R}$ is an exact functor of Waldhausen
categories.  On the other hand, since the cofiber of $\beta$ is a
model for the Eilenberg-Mac\,Lane $R$-module $H\pi_{0}R$, an
induction over Postnikov sections implies that every object in
$\aP_{R}$ is weakly equivalent to an object in $\aC_{R}^{\beta}$.
Standard arguments (e.g., \cite[VI.2.5]{EKMM} and the Whitehead
Theorem for cell $R$-modules) apply to verify that the inclusion
$\aC_{R}^{\beta}\to \aP_{R}$ satisfies the hypotheses of Waldhausen's
Approximation Theorem \cite[1.6.7]{Wald}, and thus induces a weak
equivalence of $K$-theory spectra.

The cofiber sequence of Proposition~\ref{propwaldfib} is therefore
equivalent to one of the form (cf. \cite[6.8ff]{BaasDundasRognes})
\[
K(\aP)\to K(R)\to K(R[\beta^{-1}])\to \Sigma K(\aP).
\]
Applying the D\'{e}vissage Theorem we get a weak equivalence
$K(HZ)\simeq K'(Z)\simeq K(\aP)\simeq K(\aC_{R}^{\beta})$.  To
complete the proof of the Localization Theorem we need to identify the
composite map $K(HZ)\to K(R)$ with the transfer map $K(HZ)\to K(R)$
induced from the forgetful functor.  We begin by reviewing this
transfer map.

The transfer map arises because $HZ$ is weakly equivalent to a finite
cell $R$-module; as a consequence all finite cell $HZ$-modules are
weakly equivalent to finite cell $R$-modules.  To put this into the
context of Waldhausen categories and exact functors, let $\aM_{R}$
(resp. $\aM_{HZ}$) be the full subcategory of $R$-modules
(resp. $HZ$-modules) that are weakly equivalent to finite cell
modules.  We make $\aM_{R}$ and $\aM_{HZ}$ Waldhausen categories with
w-cofibrations the Hurewicz cofibrations and weak equivalences the
usual weak equivalences.  Then the forgetful functor from $HZ$-modules
to $R$-modules restricts to an exact functor $\aM_{HZ}\to \aM_{R}$.
The inclusions $\aC_{R}\to \aM_{R}$ and $\aC_{HZ}\to \aM_{HZ}$ are
exact functors.  Again, a standard argument \cite[VI.3.5]{EKMM} with 
Waldhausen's Approximation Theorem shows that these inclusions induce 
weak equivalences of $K$-theory spectra. 
The transfer map $K(HZ)\to K(R)$ is the induced map
\[
K(HZ)\simeq K(\aM_{HZ})\to K(\aM_{R})\simeq K(R).
\]

The Waldhausen category $\aP_{R}$ is a Waldhausen subcategory of
$\aM_{R}$, and the exact functor $\aC_{R}^{\beta}\to \aM_{R}$ lands in
$\aP_{R}$ as does the exact functor $\aM_{HZ}\to \aM_{R}$.  We
therefore obtain the following commutative diagram of exact functors.
\[
\xymatrix{%
\aC_{HZ}\ar[d]&\aC^{\beta}_{R}\ar[r]\ar[d]&\aC_{R}\ar[d]\\
\aM_{HZ}\ar[r]&\aP_{R}\ar[r]&\aM_{R}
}
\]
This induces the following commutative diagram of $K$-theory spectra,
with the arrows marked ``$\simeq$'' weak equivalences.  (See also
Remark~\ref{remsmall}.) 
\[
\xymatrix{%
K(HZ)\ar[d]_{\simeq}&K(\aC_{R}^{\beta})\ar[d]_{\simeq}\ar[r]&K(R)\ar[d]^{\simeq}\\
K(\aM_{HZ})\ar[r]_{\simeq}&K(\aP_{R})\ar[r]&K(\aM_{R})
}
\]
Here the lower left map $K(\aM_{HZ})\to K(\aP_{R})$ is a weak equivalence
by the D\'{e}vissage Theorem above, since $\aM_{HZ}=\aP_{HZ}$.   This
proves the following proposition. 

\begin{prop}\label{propidentleft}
The transfer map $K(HZ)\to K(R)$ factors in the stable category as a
weak equivalence $K(HZ)\to K(\aC_{R}^{\beta})$ and the map
$K(\aC_{R}^{\beta})\to K(R)$ in Proposition~\ref{propwaldfib}. 
\end{prop}

Finally, the Localization Theorem of the introduction follows
immediately from Propositions~\ref{propwaldfib}, \ref{propidentright},
and~\ref{propidentleft}.  We close this section with a remark on smallness.

\begin{rem}\label{remsmall}
The Waldhausen categories $\aP_{R}$ and $\aM_{R}$ discussed above can
be replaced by small Waldhausen categories with equivalent
$K$-theory spectra.  For any cardinal $\kappa$, let $\aP_{\kappa}$ and
$\aM_{\kappa}$ be the full subcategories of $\aP_{R}$ and $\aM_{R}$
(respectively) consisting of those objects $M$ such that the
underlying sets of the underlying spaces of the underlying prespectrum
of $M$ are subsets of the power-set of $\kappa$.  Then $\aP_{\kappa}$
is a small category.  When $\kappa$ is bigger than the continuum and
the cardinality of the underlying sets of $R$, then $\aP_{\kappa}$ and
$\aM_{\kappa}$ contain a representative of every weak equivalence
class in $\aP_{R}$ and $\aM_{R}$, respectively.  Furthermore, $\aP_R$
and $\aM_R$ are closed under pushouts along Hurewicz cofibrations and
pullbacks along Hurewicz fibrations, and closed up to natural
isomorphism under the functors $(-)\sma X$ and $F(X,-)$ for finite
cell complexes $X$.

Let $\aC\aP_{\kappa}$ be the Waldhausen category whose objects are the
objects $M$ of $\aP_{\kappa}$ together with the structure of a cell
complex on $M$, whose morphisms are the maps of $R$-modules, whose
w-cofibrations are the maps that are isomorphic to inclusions of cell
subcomplexes, and whose weak equivalences are the usual weak
equivalences.  Let $\aC\aM_{\kappa}$ be the analogous category for
$\aM_{\kappa}$.  When $\kappa$ is bigger than the continuum and the
cardinality of the underlying sets of $R$, then for any cardinal
$\lambda \geq \kappa$, the functors $\aC\aP_{\kappa}\to \aP_{\lambda}$
and $\aC\aM_{\kappa}\to \aM_{\lambda}$ are exact and satisfy the
hypotheses of Waldhausen's Approximation Theorem \cite[1.6.7]{Wald}.
It follows that these functors induce equivalences of $K$-theory
spectra.  The reader unwilling to consider the $K$-theory of
Waldhausen categories that are not small can therefore use
$K(\aP_{\kappa})$ in place of $K(\aP_{R})$, $K(\aM_{\kappa})$ in place
of $K(\aM_{R})$, etc.
\end{rem}

\section{The $\Spdot$ Construction and the Proof of
Theorem~\ref{thmcoffib}}\label{seccoffib} 

For the Waldhausen categories $\aP_{m}^{n}$, the
condition of being a w-cofibration consists of both a point-set
requirement (HEP condition) and a homotopical requirement (injectivity
on $\pi_{n}$).  It is convenient for the arguments in the next
section to separate out these two requirements.  We do that in this
section by describing a variant $\Spdot$ of the $\Sdot$ construction
defined in terms of ``homotopy cocartesian'' squares instead of
pushouts of w-cofibrations.  For a large class of Waldhausen
categories including $\aP$, $\aP^{\op}$, and the categories
$\aP_{m}^{n}$, the $\Spdot$ construction is equivalent to the $\Sdot$
construction and can therefore be used in its place to construct
algebraic $K$-theory.  In fact, since the notions of homotopy
cocartesian and homotopy cartesian agree in $\aP$,
Theorem~\ref{thmcoffib}, which compares the algebraic $K$-theory of
$\aP$ (defined in terms of cofibrations) with that of $\aP^{\op}$
(defined in terms of fibrations), then follows as an easy consequence.

The hypothesis on a Waldhausen category we use is a weak version of
``functorial factorization''.  The Waldhausen category $\aC$ admits
\term{functorial factorization} when any map $f\colon A\to B$ in $\aC$
factors as a w-cofibration followed by a weak equivalence
\[
\xymatrix{%
A\ar@{ >->}[r]\ar@{..>}@/_1em/[rr]_{f}&Tf\ar[r]^{\sim}&B,
}
\]
functorially in $f$ in the category $\Ar\aC$ of arrows in $\aC$.  In
other words, given the map $\phi$ of arrows on the 
left (i.e., commuting diagram),
\[
\xymatrix{%
A\ar[r]^{f}\ar[d]_{a}\ar@{{}{}{}}[dr]|{\phi}&B\ar[d]^{b}
&&A\ar@{ >->}[r]\ar[d]_{a}&Tf\ar[r]^{\sim}\ar[d]_{T\phi}&B\ar[d]^{b}\\
A'\ar[r]_{f'}&B'&&A'\ar@{ >->}[r]&Tf'\ar[r]^{\sim}&B'\\
}
\]
we have a map $T\phi$ that makes the diagram on the right commute and
that satisfies the usual identity and composition relations,
$T\id_{f}=\id_{Tf}$ and $T(\phi'\circ \phi)=T\phi'\circ T\phi$.  For
example, the Waldhausen categories $\aP$ and $\aP^{\op}$ admit
functorial factorizations using the usual mapping cylinder
and mapping path-space constructions
\begin{equation}\label{eqmn}
Mf=(A\sma I_{+})\cup_{A}B\text{ for }\aP,\qquad
Nf=B^{I_{+}}\times_{B}A\text{ for }\aP^{\op}.
\end{equation}

Functorial factorization generalizes Waldhausen's notion of ``cylinder
functor satisfying the cylinder axiom'', but is still not quite
general enough to apply to the Waldhausen categories $\aP_{m}^{n}$ for
$n<\infty$, where a map is weakly equivalent to a w-cofibration only
when it is injective on $\pi_{n}$.  This leads to the following
definition.

\begin{defn}\label{defweakcof}\label{defffwc}
Let $\aC$ be a Waldhausen category.  Define a map $f\colon A\to B$
in $\aC$ to be a \term{weak w-cofibration} if is weakly equivalent in
$\Ar\aC$ (by a zigzag) to a w-cofibration.  We say that $\aC$ admits
\term{functorial factorization of weak w-cofibrations} (FFWC) when
weak w-cofibrations can be factored functorially (in $\Ar\aC$) as a
w-cofibration followed by a weak equivalence.
\end{defn}

Recall that a full subcategory $\aB$ of a Waldhausen category $\aC$ is
called a \term{Waldhausen subcategory} when it forms a Waldhausen
category with weak equivalences the weak equivalences of $\aC$ and
with w-cofibrations the w-cofibrations $A\to B$ in $\aC$ (between
objects $A$ and $B$ of $\aB$) whose quotient $B/A=B\cup_{A}*$ is in
$\aB$.  In particular, it is straight-forward to check that a full
subcategory of a Waldhausen category that is closed under extensions
is a Waldhausen subcategory; examples include the subcategories
$\aP_{m}^{n}$ of $\aP$.  We say that the Waldhausen subcategory $\aB$
is \term{closed} if every object of $\aC$ weakly equivalent to an
object of $\aB$ is an object of $\aB$.  The advantage of the
hypothesis of FFWC over the hypothesis of functorial factorization is
the following proposition.

\begin{prop}\label{propwcofpf}
If $\aB$ is a closed Waldhausen subcategory of a Waldhausen category $\aC$
that admits FFWC, then $\aB$ admits FFWC.  Moreover, a weak w-cofibration
$f \colon A \to B$ in $\aC$ between objects in $\aB$ is a weak
w-cofibration in $\aB$ if and only if the map $A \to Tf$ in
its factorization in $\aC$ is a w-cofibration in $\aB$.
\end{prop}

\begin{proof}
Let $f\colon A\to B$ be a weak w-cofibration in $\aC$ between objects
in $\aB$.  Then $f$ is weakly equivalent by a zigzag in $\aB$ to a
w-cofibration $f'\colon A'\to B'$ 
in $\aC$,
\[
\xymatrix{%
A'\ar@{ >->}[d]_{f'}\ar[r]^{\sim}&A_{1}\ar[d]_{f_{1}}&\dotsb\ar[l]_{\sim}&\ar[l]_{\sim}
A_{n}\ar[r]^{\sim}\ar[d]_{f_{n}}&A\ar[d]_{f}\\
B'\ar[r]_{\sim}&B_{1}&\dotsb\ar[l]^{\sim}
&\ar[l]^{\sim}B_{n}\ar[r]_{\sim}&B.
}
\]
Applying functorial factorization in $\aC$, we get weak equivalences
\[ 
\xymatrix@C=2em{%
B'/A'
&Tf'/A'\ar[r]^{\sim}\ar[l]_{\sim}
&Tf_{1}/A_{1}&\dotsb\ar[l]_(.4){\sim}
&\ar[l]_{\sim} Tf_{n}/A_{n}\ar[r]^{\sim}&Tf/A,
}
\]
which imply that the map $A\to Tf$ is a w-cofibration in $\aB$ if and
only if $f\colon A\to B$ is a weak w-cofibration in $\aB$ .
\end{proof}

In general, the weak w-cofibrations do not necessarily form a
subcategory of $\aC$; the composition of two weak w-cofibrations might
not be a weak w-cofibration.  However, in the presence of FFWC the
weak w-cofibrations are well-behaved.

\begin{prop}
Let $\aC$ be a Waldhausen category that admits FFWC.  If $f \colon A \to B$
and $g \colon B \to C$ are weak w-cofibrations in $\aC$, then $g \circ
f \colon  A \to C$ is a weak w-cofibration in $\aC$.
\end{prop}

\begin{proof}
Applying functorial factorization in $\aC$, we can factor $f$ and $g$
as the composites
\[
\xymatrix{%
A\ar@{ >->}[r]&Tf\ar[r]^{\sim}&B &&
B\ar@{ >->}[r]&Tg\ar[r]^{\sim}&C.
}
\]
The composite map $h \colon Tf \to B \to Tg$ is a weak
w-cofibration, as it is weakly equivalent in $\Ar\aC$ to $g$.  We can
therefore apply the factorization functor to $h$ to obtain a
w-cofibration $Tf\to Th$ and weak equivalence $Th\to Tg$.
The composite w-cofibration $A\to Tf \to Th$ is
then weakly equivalent to $g \circ f$.
\end{proof}

A square diagram in a Waldhausen category (for example) is called
\term{homotopy cocartesian} if it is weakly equivalent (by a zigzag)
to a pushout square where one of the parallel sets of arrows 
consists of w-cofibrations.  The corresponding set of parallel arrows in the
original square then consists of weak w-cofibrations.  When a
Waldhausen category admits FFWC, there is a good criterion in terms of
the factorization functor $T$ for detecting homotopy cocartesian
squares.  We state it for a Waldhausen category that is
\term{saturated}, i.e., one whose weak equivalences satisfy the
``two-out-of-three property'' (see, for example, \cite[1.2]{Wald}). 
The proof is similar to that of Proposition~\ref{propwcofpf}.

\begin{prop}\label{prophococart}
Let $\aC$ be a saturated Waldhausen category admitting FFWC with functor
$T$.  A commutative diagram 
\[
\xymatrix@-1pc{%
A\ar[r]^{f}\ar[d]&B\ar[d]\\
C\ar[r]&D
}
\]
with $f$ a weak w-cofibration is homotopy cocartesian if and only
if the map $Tf\cup_{A}C\to D$ is a weak equivalence.
\end{prop}

Likewise, for saturated Waldhausen categories admitting FFWC, 
homotopy cocartesian squares have many of the usual expected
properties.  We summarize the ones we need in the following
proposition; its proof is a straight-forward application of the
previous proposition.

\begin{prop}\label{proppush}
Let $\aC$ be a saturated Waldhausen category admitting FFWC.
\begin{enumerate}
\item Given a commutative cube
\[
\xymatrix@-1.5pc{%
A'\ar[rr]\ar[dd]\ar[dr]&&B'\ar[dd]\ar[dr]\\
\relax&A\ar'[r][rr]\ar'[d][dd]&&B\ar[dd]\\
C'\ar[rr]\ar[dr]&&D'\ar[dr]\\
\relax&C\ar[rr]&&D
}
\]
with the $(A,B,C,D)$-face and $(A',B',C',D')$-face homotopy
cocartesian, if the maps $A'\to A$, 
$B'\to B$, and $C'\to C$ weak equivalences, 
then the map $D'\to D$ is a weak equivalence.
\item Given a commutative diagram
\[
\xymatrix@-1pc{%
A\ar[r]\ar[d]&B\ar[d]\ar[r]&X\ar[d]\\ 
C\ar[r]&D\ar[r] &Y \\
}
\]
with the square $(A,B,C,D)$ homotopy cocartesian, if either $A\to C$
is a weak w-cofibration or both $A\to B$ and $B\to X$ are weak
w-cofibrations, then the $(A,X,C,Y)$
square is homotopy cocartesian if and only if the $(B,X,D,Y)$ square
is homotopy cocartesian.
\end{enumerate}
\end{prop}

We use the concept of weak w-cofibration and the previous propositions
to build a homotopical variant of the $\Sdot$ construction.  First, we recall
the $\Sdot$ construction in detail.  Waldhausen's $\Sdot$ construction
produces a simplicial Waldhausen category $\Sdot\aC$ from a Waldhausen
category $\aC$ and is defined as follows.  Let $\Ar[n]$ denote the
category with objects $(i,j)$ for $0\leq i\leq j\leq n$ and a unique
map $(i,j)\to (i',j')$ for $i\leq i'$ and $j\leq j'$.  $\Sdot[n]\aC$
is defined to be the full subcategory of the category of functors
$A\colon \Ar[n]\to \aC$ such that:
\begin{itemize}
\item $A_{i,i}=*$ for all $i$, 
\item The map $A_{i,j}\to A_{i,k}$ is a w-cofibration for all $i \leq
j \leq k$, and
\item The diagram
\[  \xymatrix@-1pc{%
A_{i,j}\ar[r]\ar[d]&A_{i,k}\ar[d]\\A_{j,j}\ar[r]&A_{j,k}
} \]
is a pushout square for all $i \leq j \leq k$, 
\end{itemize}
where we write $A_{i,j}$ for $A(i,j)$.  The last two conditions can be
simplified to the hypothesis that each map $A_{0,j}\to A_{0,j+1}$ is a
w-cofibration and the induced maps $A_{0,j}/A_{0,i}\to A_{i,j}$ are
isomorphisms.  This becomes a Waldhausen category by defining a map
$A\to B$ to be a weak equivalence when each $A_{i,j}\to B_{i,j}$ is a
weak equivalence in $\aC$, and to be a w-cofibration when each
$A_{i,j}\to B_{i,j}$ and each induced map
$A_{i,k}\cup_{A_{i,j}}B_{i,j}\to B_{i,k}$ is a w-cofibration in $\aC$.

\begin{defn}
Let $\aC$ be a saturated Waldhausen category that admits FFWC.  Define
$\Spdot[n]\aC$ 
to be the full subcategory of functors $A\colon \Ar[n]\to \aC$ such that:
\begin{itemize}
\item The initial map $*\to A_{i,i}$ is a weak equivalence for all $i$, 
\item The map $A_{i,j}\to A_{i,k}$ is a weak w-cofibration for all $i
\leq j \leq k$, and 
\item The diagram
\[  \xymatrix@-1pc{%
A_{i,j}\ar[r]\ar[d]&A_{i,k}\ar[d]\\A_{j,j}\ar[r]&A_{j,k}
} \]
is a homotopy cocartesian square for all $i \leq j \leq k$.
\end{itemize}
We define a map $A\to B$ to be a
weak equivalence when each $A_{i,j}\to B_{i,j}$ is a weak equivalence
in $\aC$, and to be a w-cofibration when each $A_{i,j}\to B_{i,j}$ is
a w-cofibration in $\aC$ and each induced map
$A_{i,k}\cup_{A_{i,j}}B_{i,j}\to B_{i,k}$ is a weak w-cofibration in
$\aC$.
\end{defn}

The following theorem in particular allows us to iterate the
$\Spdot$ construction. 

\begin{thm}\label{thmspdotinherit}
Let $\aC$ be a saturated Waldhausen category that admits FFWC; then
$\Spdot[n]\aC$ is also a saturated Waldhausen category that admits FFWC.
\end{thm}

\begin{proof}

The proof that $\Spdot[n]\aC$ is a saturated Waldhausen category
follows the same outline as the proof that $\Sdot[n]\aC$ is a
saturated Waldhausen category, substituting
Propositions~\ref{prophococart} and~\ref{proppush} for the homotopy
cocartesian squares in place of the usual arguments for pushout
squares. To show
that $\Spdot[n]\aC$ admits FFWC, we construct a factorization functor
$T$ as follows. Given $f\colon A\to 
B$ a weak w-cofibration in $\Spdot[n] \aC$, let
$(Tf)_{i,j}=T(f_{i,j})$.  This then factors $f$ as a map $A\to Tf$
followed by a weak equivalence $Tf\to B$; we need to check that the
map $A\to Tf$ is a w-cofibration.  Since by construction the maps
$A_{i,j}\to Tf_{i,j}$ are w-cofibrations, we just need to check that
the maps $A_{i,k}\cup_{A_{i,j}}Tf_{i,j}\to Tf_{i,k}$ are weak
w-cofibrations. Since by hypothesis $f$ is a weak w-cofibration in
$\Spdot[n] \aC$, it is weakly equivalent (by a zigzag) to a
w-cofibration $f'\colon A'\to B'$ in $\Spdot[n] \aC$.  It follows that
the maps
\[
A_{i,k}\cup_{A_{i,j}}Tf_{i,j}\to Tf_{i,k} \qquad\text{and}\qquad
A'_{i,k}\cup_{A'_{i,j}}Tf'_{i,j}\to Tf'_{i,k}
\]
are weakly equivalent (by the corresponding zigzag).  Since
$A'_{i,j}\to B'_{i,j}$ is a w-cofibration, the canonical map
\[
\xymatrix@R-1em{%
A'_{i,k}\cup_{A'_{i,j}}Tf'_{i,j}\ar[d]\ar[rr]^{\sim}
&&A'_{i,k}\cup_{A'_{i,j}}B'_{i,j}\ar[d]\\
Tf'_{i,k}\ar[rr]_{\sim}&&B'_{i,k}
}
\]
is a weak equivalence.  The map $A'_{i,k}\cup_{A'_{i,j}}B'_{i,j}\to
B'_{i,k}$ is a weak w-cofibration by hypothesis, and so the maps
$A'_{i,k}\cup_{A'_{i,j}}Tf'_{i,j}\to Tf'_{i,k}$ and
$A_{i,k}\cup_{A_{i,j}}Tf_{i,j}\to Tf_{i,k}$ are therefore weak
w-cofibrations. 
\end{proof}

Both $\Sdot\aC$ and $\Spdot\aC$ become simplicial Waldhausen
categories with face map $\partial_{i}$ deleting the $i$-th row and
column, and degeneracy map $s_{i}$ repeating the $i$-th row and
column.  For each $n$, we denote the nerve of the category of weak
equivalences in $\Sdot[n]\aC$ by $w\subdot\Sdot[n]\aC$.  As $n$
varies, $w\subdot \Sdot\aC$ assembles into a bisimplicial set, and we
denote the geometric realization by $|w\subdot\Sdot\aC|$.  By
definition, the algebraic $K$-theory space of $\aC$ is $\Omega
|w\subdot\Sdot\aC|$.  The algebraic $K$-theory spectrum is obtained by
iterating the $\Sdot$ construction: Its $n$-th space is $|w\subdot
\Sdot^{(n)}\aC|$.  The following theorem therefore implies that
algebraic $K$-theory can be constructed from the $\Spdot$
construction.

\begin{thm}\label{thmspdot}
Let $\aC$ be a saturated Waldhausen category that admits FFWC.  The
inclusion of $\Sdot\aC$ in $\Spdot\aC$ 
is an exact simplicial functor and the induced map of bisimplicial
sets $w\subdot\Sdot\aC \to w\subdot\Spdot\aC$ is a weak equivalence. 
\end{thm}

\begin{proof}
The fact that $\Sdot\aC\to \Spdot\aC$ is an exact simplicial functor
is clear from the definition of the simplicial structure and the weak
equivalences and w-cofibrations in each category.  To see that it
induces a weak equivalence on nerves, it suffices to show that for
each $n$, the map $w\subdot \Sdot[n]\aC\to w\subdot \Spdot[n]\aC$ is a
weak equivalence.  We do this by constructing a functor $\Phi\colon
\Spdot[n]\aC\to\Sdot[n]\aC$ and natural weak equivalences relating
$\Phi$ and the identity in both $\Spdot[n]\aC$ and $\Sdot[n]\aC$.
In $\Spdot[n]\aC$, we construct a zigzag of natural weak equivalences of
the form 
\[
\xymatrix{%
\Phi &\Theta \ar[l]_{q}^{\simeq}\ar[r]^{\epsilon}_{\simeq}&\Id
}
\]
and in $\Sdot[n]\aC$, we construct a weak equivalence $\Phi \to \Id$.

We define the functor $\Theta$ inductively as follows.
Given an object $A$ in $\Spdot[n]\aC$, let $\Theta_{0,0}A=A_{0,0}$ and let
$\epsilon_{0,0}\colon \Theta_{0,0}A\to A_{0,0}$ be the identity map.
Given $\Theta_{0,j}A$ and $\epsilon_{0,j}\colon \Theta_{0,j}A\to
A_{0,j}$, define $\Theta_{0,j+1}A$ to be $Tf$ and $\epsilon_{0,j+1}\colon
Tf\to A_{0,j+1}$ to be the canonical weak equivalence, where $f\colon
\Theta_{0,j}A\to A_{0,j+1}$ is the composite of $\epsilon_{0,j}$ with the
map $A_{0,j}\to A_{0,j+1}$ in the structure of the functor $A$.  As a
consequence of the construction, we have w-cofibrations $\Theta_{0,j}A\to
\Theta_{0,j+1}A$.  Next we construct the diagonal objects $\Theta_{i,i}A$;
given $\Theta_{j,j}A$ and $\epsilon_{j,j}\colon \Theta_{j,j}A\to
A_{j,j}$, define $\Theta_{j+1,j+1}A$ to be $Tf$ and $\epsilon_{j+1,j+1}\colon
Tf\to A_{0,j+1}$ to be the canonical weak equivalence, for $f\colon
\Theta_{j,j}A\to A_{j,j}\to A_{j+1,j+1}$.  Finally, for $0<i<j$, let
\[
\Theta_{i,j} A=\Theta_{0,j} A\cup_{\Theta_{0,i}A} \Theta_{i,i} A
\]
and let $\epsilon_{i,j}\colon \Theta_{i,j} A\to A_{i,j}$ be the map
induced by the universal property of the pushout.  By
Proposition~\ref{proppush}, the maps $\epsilon_{i,j}$ are weak
equivalences.  Thus, $\Theta$ defines a functor from $\Spdot[n]\aC$ to
itself, and $\epsilon$ defines a natural transformation from $\Theta$
to the identity.

We define $\Phi$ and $q$ by setting
$\Phi_{i,i} A = *$,
\[
\Phi_{i,j} A = \Theta_{i,j} A / \Theta A_{i,i},
\]
and $q_{i,j}\colon \Theta_{i,j}A\to \Phi_{i,j}A$ the quotient map.
Since the initial map $*\to \Theta_{i,i}A$ is a weak equivalence, so
is the final map $\Theta_{i,i}A\to *$, and thus each $q_{i,j}$ is a weak
equivalence.  The construction of $\Theta_{i,j} A$ 
as a pushout for $0<i<j$ implies that the induced map from the quotient
$\Phi_{0,j}A/\Phi_{0,i}A$ to $\Phi_{i,j}A$ is an isomorphism.  It
follows that $\Phi$ 
defines a functor from $\Spdot[n]\aC$ to $\Sdot[n]\aC$ and $q$ a natural
weak equivalence in $\Spdot[n]\aC$ from $\Theta$ to $\Phi$.

Finally, if $A$ is an object of $\Sdot[n]\aC$, then $A_{i,i}=*$ for all
$i$, and the map $\epsilon \colon \Theta A\to A$ factors through the
quotient $\Phi A$.  This defines a natural weak equivalence from
$\Phi$ to the identity in $\Sdot[n]\aC$.
\end{proof}

Theorem~\ref{thmcoffib} is an easy consequence of
Theorem~\ref{thmspdot}: Because a square in $\aP$ is homotopy
cocartesian if and only if it is homotopy cartesian (homotopy
cocartesian in $\aP^{\op}$), the categories $\Spdot[n]\aP$ and
$\Spdot[n]\aP^{\op}$ are (contravariantly) isomorphic for each $n$ by
the isomorphism that renumbers the diagram by the involution
$(i,j)\mapsto (n-j,n-i)$.  The bisimplicial sets $w\subdot \Spdot \aP$
and $w\subdot \Spdot\aP^{\op}$ then only differ in the orientation of
the simplices, and the spaces $\Omega |w\subdot \Spdot \aP|$ and
$\Omega |w\subdot \Spdot \aP^{\op}|$ are homeomorphic.  Comparing
iterates of the $\Spdot$ construction and the suspension maps, we
obtain an equivalence of $K$-theory spectra.

\section{A Reduction of Theorem~\ref{thmbigthm}}\label{secbigthm}

The purpose of this section is to reduce Theorem~\ref{thmbigthm} to
Theorem~\ref{thmemreduc} below.  That theorem is similar in spirit to
Theorem~\ref{thmpizero}, and both are proved with closely related
arguments in the next section.

We begin by reducing to the subcategory of connective objects.

\begin{lem}\label{lemconn}
The inclusion $\aP_{0}^{\infty}\to \aP$ induces an equivalence of $K$-theory.
\end{lem}

\begin{proof}
Since $\aP=\colim \aP_{-m}^{\infty}$ and $\Sdot \aP=\colim \Sdot
\aP_{-m}^{\infty}$, it suffices to show that the inclusions
$\aP^{\infty }_{-m}\to \aP^{\infty}_{-m-1}$ induce equivalences of
$K$-theory.  Suspension gives an exact functor $\aP_{-m-1}^{\infty}\to
\aP_{-m}^{\infty}$. The composite endomorphisms on $\aP_{-m}^{\infty}$
and $\aP_{-m-1}^{\infty}$ are the suspension functors, which on the
$K$-theory spectra induce multiplication by $-1$ \cite[1.3.2.4]{Wald} and
in particular are weak equivalences.
\end{proof}

We also have $\aP_{0}^{\infty}=\colim \aP_{0}^{n}$ and $\Sdot
\aP_{0}^{\infty} =\colim \Sdot \aP_{0}^{n}$, which gives the following proposition.

\begin{prop}\label{proptel}
The spectrum $K\aP_{0}^{\infty}$ is equivalent to the telescope of the
sequence of maps $K\aP_{0}^{0}\to \dotsb \to K\aP_{0}^{n}\to
K\aP_{0}^{n+1}\to \dotsb$.
\end{prop}

The proof of Theorem~\ref{thmbigthm} will then be completed by showing
that the maps $\aP_{0}^{n}\to \aP_{0}^{n+1}$ induce weak equivalences of
$K$-theory spectra, which we do by studying the cofiber.  According to
Waldhausen \cite[1.5.6]{Wald}, for a Waldhausen category $\aC$ and
a Waldhausen subcategory $\aB$, the spectrum-level cofiber of $K\aB\to
K\aC$ is the $K$-theory spectrum of the simplicial Waldhausen category
$F\subdot(\aC,\aB)$, defined as follows.  The Waldhausen category
$F_{q}(\aC,\aB)$ has as objects the sequences of $q$ composable
w-cofibrations in $\aC$
\[
\xymatrix@C=1.5em{%
C_{0}\ar@{ >->}[r]&\dotsb\ar@{ >->}[r]&C_{i}\ar@{ >->}[r]&\dotsb\ar@{ >->}[r]&C_{q}
}
\]
such that the quotients $C_{j}/C_{i}$ are in $\aB$ for all $0\leq i<j\leq q$ and
has as morphisms commutative diagrams.  A map $C\to D$ is a 
weak equivalences when the maps $C_{i}\to D_{i}$ are weak equivalences
in $\aC$ for all $i$, and is a w-cofibration when the maps $C_{i}\to D_{i}$ are
w-cofibrations for all $i$ and the maps $C_{j}\cup_{C_{i}}D_{i}\to
D_{j}$ are w-cofibrations for all $i<j$.  The face map $\partial_{i}$
deletes the object $C_{i}$ and uses the composite w-cofibration
$C_{i-1}\to C_{i+1}$, and the degeneracy map $s_{i}$ repeats the object
$C_{i}$, inserting the identity map.  We use the following variant of
this construction.

\begin{defn}
For a saturated Waldhausen category $\aC$ that admits FFWC and a
closed Waldhausen subcategory $\aB$, let
$F'_{q}(\aC,\aB)$ denote the category that has as objects the
sequences of $q$ composable weak w-cofibrations in $\aC$ 
\[
\xymatrix@C=1.5em{%
C_{0}\ar[r]&\dotsb\ar[r]&C_{i}\ar[r]&\dotsb\ar[r]&C_{q}
}
\]
such that  for all $0\leq i<j\leq q$ there exists an object $B_{i,j}$ in
$\aB$ and a map $C_{j}\to B_{i,j}$ for which the square
\[  \xymatrix@-1pc{%
C_{i}\ar[r]\ar[d]&C_{j}\ar[d]\\
{*}\ar[r]&B_{i,j}
} \]
commutes and is homotopy cocartesian.  A morphism $C\to D$ in 
$F'_{q}(\aC,\aB)$ is a  commutative diagram
\[
\xymatrix@=1.5em{%
C_{0}\ar[r]\ar[d]&\dotsb\ar[r]&C_{i}\ar[d]\ar[r]&\dotsb\ar[r]&C_{q}\ar[d]\\
D_{0}\ar[r]&\dotsb\ar[r]&D_{i}\ar[r]&\dotsb\ar[r]&D_{q}
}
\]
We make
$F'_{q}(\aC,\aB)$ a Waldhausen category  by declaring a map $C\to D$
to be a weak equivalence when the maps $C_{i}\to D_{i}$ are weak equivalences
in $\aC$ for all $i$, and to be a w-cofibration when the maps $C_{i}\to D_{i}$ are
w-cofibrations for all $i$ and the maps $C_{j}\cup_{C_{i}}D_{i}\to
D_{j}$ are weak w-cofibrations for all $i<j$.  
\end{defn}

Note that in the above definition, the choice of ``$B_{i,j}$'' in
$\aB$ is not part of the structure of the object $C$.  The hypothesis
involving the existence of the objects $B_{i,j}$ is a precise way of
making sense of the condition that the ``homotopy cofiber'' of
$C_{i}\to C_{j}$ lies in $\aB$ for an arbitrary Waldhausen category
$\aC$.
We assemble
$F'\subdot(\aC,\aB)$ into a simplicial Waldhausen category just as
$F\subdot(\aC,\aB)$ above: The face map $\partial_{i}$
deletes the object $C_{i}$ and uses the composite weak w-cofibration
$C_{i-1}\to C_{i+1}$; the degeneracy map $s_{i}$ repeats the object
$C_{i}$, inserting the identity map.  The inclusion $F\subdot(\aC,\aB)\to
F'\subdot(\aC,\aB)$ is an exact simplicial functor.

We have stated the definitions of $F$ and $F'$ for general $\aC,\aB$
rather than just for $\aP_{0}^{n+1},\aP_{0}^{n}$ in order to apply a
standard trick of Waldhausen \cite{Wald} of commuting these
constructions past another construction, in this case the $\Spdot$
construction.  When $\aC$ is a saturated Waldhausen category that
admits FFWC, the 
Waldhausen category $\Spdot[p]F_{q}(\aC,\aB)$ is isomorphic to
$F_{q}(\Spdot[p]\aC,\Spdot[p]\aB)$ and the Waldhausen category
$\Spdot[p]F'_{q}(\aC,\aB)$ is isomorphic to
$F'_{q}(\Spdot[p]\aC,\Spdot[p]\aB)$.  The argument for
Theorem~\ref{thmspdot} also shows that the inclusion $w\subdot
F_{q}(\aC,\aB)\to w\subdot F'_{q}(\aC,\aB)$ is a homotopy equivalence
of simplicial sets for all $q$.  Since $\Spdot[p]\aC$ is also a
saturated Waldhausen category that admits FFWC
and $\Spdot[p]\aB$ is a closed Waldhausen subcategory of $\Spdot[p]\aC$, it
follows that the inclusion $w\subdot
F_{q}(\Spdot[p]\aC,\Spdot[p]\aB)\to w\subdot
F'_{q}(\Spdot[p]\aC,\Spdot[p]\aB)$ is a homotopy equivalence of
simplicial sets for all $p,q$.  This proves the following proposition.

\begin{prop}
Let $\aC$ be a saturated Waldhausen category that admits FFWC and $\aB$ a
closed Waldhausen subcategory. 
Then the inclusion $F\subdot(\aC,\aB)\to
F'\subdot(\aC,\aB)$ induces a weak equivalence of $K$-theory spectra. Thus,
we have a cofibration sequence
\[
K\aB\to K\aC \to KF'\subdot(\aC,\aB) \to \Sigma K\aB
\]
in the stable category.
\end{prop}

The remainder of the proof of Theorem~\ref{thmbigthm} is to show that
$w\subdot\Spdot F'\subdot(\aP_{0}^{n+1},\aP_{0}^{n})$ is weakly contractible
for all $n \geq 0$.  The first step is to eliminate one of the
simplicial directions.

\begin{defn}
Let $\cU$ be the subcategory of $\aP_{0}^{n+1}$ consisting of all
objects and those maps that induce an isomorphism on $\pi_{n+1}$ and
an injection in $\pi_{n}$.  Let $u\Spdot[p]\aP_{0}^{n+1}$ denote the
subcategory of $\Spdot[p]\aP_{0}^{n+1}$ consisting of all objects and
those maps $f$ in $\Spdot[p]\aP_{0}^{n+1}$ such that $f_{i,j}$ is in
$\cU$ for all $0\leq i < j\leq p$.
\end{defn}

To avoid a possible point of confusion, note that $\cU$ and
$u\Spdot[p]\aP_{0}^{n+1}$ are not categories of weak equivalences in
the sense of Waldhausen.  Rather, the point of the category $\cU$ is
that a map $f\colon A\to B$ in 
$\aP_{0}^{n+1}$ has homotopy cofiber in $\aP_{0}^{n}$ if and only $f$
is a map in $\cU$.  Writing $u\subdot \Spdot[p]\aP_{0}^{n+1}$ for the
nerve of $u\Spdot[p]\aP_{0}^{n+1}$, 
the set of $q$-simplexes $u_{q} \Spdot[p]\aP_{0}^{n+1}$ is precisely
the set of objects of $F_{q}'(\Spdot[p]\aP_{0}^{n+1},\Spdot[p]\aP_{0}^{n})$.
The following lemma is a special case of Waldhausen's ``Swallowing
Lemma'' \cite[1.6.5]{Wald}. 

\begin{lem}
The inclusion 
\[
u\subdot \Spdot[p] \aP_{0}^{n+1}=w_{0}F'\subdot
(\Spdot[p] \aP_{0}^{n+1},\Spdot[p]\aP_{0}^{n}) \to 
w\subdot F'\subdot(\Spdot[p] \aP_{0}^{n+1},\Spdot[p]\aP_{0}^{n})
\]
is a weak equivalence for all $p$.
\end{lem}

\begin{proof}
Let $uw_{r}\Spdot[p]$ denote the category whose objects are the sequences of $r$
composable weak equivalences between objects in
$\Spdot[p]\aP_{0}^{n+1}$,
\[
\xymatrix@=1.5em{%
A_{0}\ar[r]^{\sim}&\dotsb\ar[r]^{\sim}&A_{i}\ar[r]^{\sim}&\dotsb\ar[r]^{\sim}&A_{r},
}
\]
and whose maps $A\to B$ are the commutative diagrams
\[
\xymatrix@=1.5em{%
A_{0}\ar[r]^{\sim}\ar[d]_{u}&\dotsb\ar[r]^{\sim}&A_{i}\ar[r]^{\sim}\ar[d]_{u}
&\dotsb\ar[r]^{\sim}&A_{r}\ar[d]_{u}\\
B_{0}\ar[r]_{\sim}&\dotsb\ar[r]_{\sim}&B_{i}\ar[r]_{\sim}&\dotsb\ar[r]_{\sim}&B_{r}
}
\]
with each map $A_{i}\to B_{i}$ in $u\Spdot[p]\aP_{0}^{n+1}$.  Then
for each $r$, $w_{r} F'\subdot(\Spdot[p] \aP_{0}^{n+1},\Spdot[p]\aP_{0}^{n})$ is
the nerve $u\subdot w_{r}\Spdot[p]$ of the category  $u\subdot
w_{r}\Spdot[p]$, and the map in question is the 
inclusion of $u\subdot w_{0}\Spdot[p]$ in $u\subdot
w\subdot\Spdot[p]$.  Thus, to show that it is a weak equivalence, it
suffices to show that each iterated degeneracy $u\subdot
w_{0}\Spdot[p]\to u\subdot w_{r}\Spdot[p]$ is a weak equivalence.  The
map $u\subdot w_{0}\Spdot[p]\to u\subdot w_{r}\Spdot[p]$ is the nerve
of the functor that sends $A$ to the sequence of identity maps
$A=\dotsb=A$.  This is left adjoint to the functor $uw_{r}\Spdot[p]\to
uw_{0}\Spdot[p]$ that sends a sequence $A_{0}\to\dotsb \to A_{r}$ to
its zeroth object $A_{0}$.  It follows that $u\subdot
w_{0}\Spdot[p]\to u\subdot w_{r}\Spdot[p]$ is a homotopy equivalence.
\end{proof}

At this point we could state Theorem~\ref{thmemreduc}, although it
would appear somewhat cryptic.  To explain the form it takes, consider
the further simplification of replacing the upper triangular diagrams
of $\Spdot[p]$ with a sequence of composable maps (at the cost of
losing the simplicial structure in $p$).  For a Waldhausen category
$\aC$, let $F_{p-1} \aC=F_{p-1}(\aC,\aC)$ be the Waldhausen category
whose objects are the sequences of $p-1$ composable w-cofibrations in
$\aC$.  Then the functor $w\Sdot[p]\aC\to wF_{p-1}\aC$ that sends
$(A_{i,j})$ to $A_{0,1}\to \dotsb \to A_{0,p}$ is an equivalence of
categories and therefore induces a weak equivalence of nerves.  When
$\aC$ is a saturated Waldhausen category that admits FFWC, we have the
Waldhausen category 
$F'_{p-1}\aC=F'_{p-1}(\aC,\aC)$ and the analogous functor
$w\Spdot[p]\aC\to wF'_{p-1}\aC$ induces a weak equivalence of nerves.
This argument does not apply to the map
\[
u\subdot \Spdot[p]\aP_{0}^{n+1}\to 
u\subdot F_{p-1}'\aP_{0}^{n+1},
\]
however.  The reason is that when $A_{1}\to A_{2}$ and $B_{1}\to
B_{2}$ are weak w-cofibrations in $\aP_{0}^{n+1}$ (maps that are
injective on $\pi_{n+1}$), and
\[
\xymatrix@=1.5em{%
A_{1}\ar[r]\ar[d]_{u}&A_{2}\ar[d]^{u}\\
B_{1}\ar[r]&B_{2}
}
\]
is a commutative diagram with the vertical maps in $\cU$,  the induced
map on homotopy pushouts  (with notation as in \eqref{eqmn}) 
\[
M(A_{1}\rightarrow A_{2})\longrightarrow M(B_{1}\rightarrow B_{2})
\]
is generally not a map in $\cU$.  The key is to look instead at the
``weak fibration'' analogue: If we assume instead that $A_{1}\to
A_{2}$ and $B_{1}\to B_{2}$ are maps in $\aP_{0}^{n+1}$ that induce
epimorphisms on $\pi_{0}$, then for a commutative diagram as above,
the induced map of homotopy pullbacks
\[
N(A_{1}\rightarrow A_{2})\longrightarrow N(B_{1}\rightarrow B_{2})
\]
is a map in $\cU$ by the Five Lemma.  This suggests the following definition.

\begin{defn}
Let $F^{f}_{q}\aP_{0}^{n+1}$ be the category where an object is a
sequence of $q$ composable maps $A_{0}\to \dotsb \to A_{q}$ in
$\aP_{0}^{n+1}$ that induce epimorphisms on $\pi_{0}$ and a morphism
is a commutative diagram. Let $uF^{f}_{q}\aP_{0}^{n+1}$ denote the
subcategory containing all objects of $F^{f}_{q}\aP_{0}^{n+1}$, but
consisting of only those morphisms $A\to B$ such that $A_{i}\to B_{i}$
is in $\cU$ for all $i$.
\end{defn}

For any object $A=(A_{i,j})$ in $\Spdot[p]\aP_{0}^{n+1}$, the maps
$A_{i,p}\to A_{j,p}$ induce epimorphisms on $\pi_{0}$ for all $0\leq
i<j\leq p-1$, and so we get a functor $\phi \colon
\Spdot[p]\aP_{0}^{n+1}\to F^{f}_{p-1}\aP_{0}^{n+1}$ that sends $A$ to
$A_{0,p}\to A_{1,p}\to \dotsb \to A_{p-1,p}$.  Moreover, since the
squares in $A$ are homotopy (co)cartesian, we can recover $A_{i,j}$
from $\phi A$ up to weak equivalence as the homotopy fiber of the map
$A_{i,p}\to A_{j,p}$ for all $0\leq i<j\leq p-1$.  Likewise, for any
object $(A_{i})$ in $F^{f}_{p-1}\aP_{0}^{n+1}$, the homotopy pullback
of each map
$A_{i}\to A_{j}$ is an object of $\aP_{0}^{n+1}$, and the mapping
path-space construction analogous to the construction in the proof of
Theorem~\ref{thmspdot} constructs a functor $\Phi \colon
F^{f}_{p-1}\aP_{0}^{n+1}\to \Spdot[p]\aP_{0}^{n+1}$ and natural weak
equivalences relating $\phi \circ \Phi$ and $\Phi \circ \phi$ to the
identity functors in $F^{f}_{p-1}\aP_{0}^{n+1}$ and $\Spdot[p]\aP_{0}^{n+1}$.
As mentioned above, the Five Lemma implies that a morphism $A\to B$ in
$\Spdot[p]\aP_{0}^{n+1}$ is in $u\Spdot[p]\aP_{0}^{n+1}$ if and only
if $\phi A\to \phi B$ is in $uF^{f}_{p-1}\aP_{0}^{n+1}$; thus, $\phi$
and $\Phi$ restrict to functors $\phi \colon
u\Spdot[p]\aP_{0}^{n+1}\to uF^{f}_{p-1}\aP_{0}^{n+1}$ and $\Phi \colon
uF^{f}_{p-1}\aP_{0}^{n+1} \to u\Spdot[p]\aP_{0}^{n+1}$.  As a
consequence, we obtain the following proposition.

\begin{prop}\label{propuff}
The functor $\phi \colon u\Spdot[p]\aP_{0}^{n+1}\to
uF^{f}_{p-1}\aP_{0}^{n+1}$ induces a weak equivalence of nerves for
each $p$.
\end{prop}

For any object $A$ of $\aP_{0}^{n+1}$, the $n$-connected cover
$A\ang[n]$ is an Eilenberg--Mac\,Lane $R$-module $K(\pi_{n+1}A,n+1)$,
and the map $A\ang[n]\to A$ is in $\cU$.  This suggests that the nerve
$u\subdot F^{f}_{p-1}\aP_{0}^{n+1}$ is weakly equivalent to the nerve
of the full subcategory of sequences of maps of $K(\pi,n+1)$'s (see
Lemma~\ref{lemapplymainlem} below for a precise statement).  This
further suggests that we can understand $u\subdot
F^{f}_{p-1}\aP_{0}^{n+1}$ and therefore
$u\subdot\Spdot[p]\aP_{0}^{n+1}$ in terms of $\pi_{n+1}$.  This leads
to the following definition.

\begin{defn}
Let $Z=\pi_{0}R$.  Let $M_{p}Z$ be the category whose objects are
sequences of $p-1$ composable maps of finitely generated $Z$-modules
$X_{0}\to \dotsb \to X_{p-1}$ and whose morphisms are commutative
diagrams.  Let $uM_{p}Z$ be the subcategory of $M_{p}Z$ consisting of
all objects but only those maps $X\to Y$ that are isomorphisms
$X_{i}\to Y_{i}$ for all $0\leq i\leq p-1$.
\end{defn}

We understand $uM_{0}Z$ to be the trivial category consisting of a
single object (the empty sequence of maps) with only the identity map.
We then make $uM\subdot Z$ a simplicial category as follows: For
$0\leq i\leq p-1$, the face map $\partial_{i}\colon M_{p}Z\to
M_{p-1}Z$ is obtained by dropping $X_{i}$ (and composing) and the
degeneracy map $s_{i}\colon M_{p-1}Z\to M_{p}Z$ is obtained by
repeating $X_{i}$ (with the identity map).  The face map
$\partial_{p}\colon M_{p}Z\to M_{p-1}Z$ sends $X_{0}\to \dotsb \to
X_{p-1}$ to $K_{0}\to \dotsb \to K_{p-2}$, where $K_{i}$ is the kernel
of the composite map $X_{i}\to X_{p-1}$.  The last degeneracy
$s_{p-1}\colon M_{p-1}Z\to M_{p}Z$ puts $0$ in as the last object in
the sequence.  With this definition, the functors $\psi_{p} \colon
u\Spdot[p]\aP_{0}^{n+1}\to uM_{p}Z$ that take $A=(A_{i,j})$ in
$\Spdot[p]\aP_{0}^{n+1}$ to
\[
\pi_{n+1}(A_{0,p}) \to \pi_{n+1}(A_{1,p}) \to \dotsb \to \pi_{n+1}(A_{p-1,p}) 
\]
in $M_{p}Z$ assemble to a simplicial functor $\psi \colon
u\Spdot\aP_{0}^{n+1}\to uM\subdot Z$.  Note also that for each $p$,
the functor $\psi_{p}$ factors as the composite $\pi_{n+1}\circ
\phi$.  We will exploit this factorization to prove the following
theorem in the next section.

\begin{thm}\label{thmemreduc}
The simplicial functor $\psi \colon u\Spdot\aP_{0}^{n+1}\to
uM\subdot Z$ induces a weak equivalence of nerves.
\end{thm}

Finally, to reduce Theorem~\ref{thmbigthm} to
Theorem~\ref{thmemreduc}, we need to show that the nerve $u\subdot
M\subdot Z$ is weakly contractible.  For this, consider the functors
$c_{p} \colon uM_{p}Z\to uM_{p+1}Z$ that take the sequence
$X_{0}\to \dotsb \to X_{p-1}$ to the sequence $0\to X_{0}\to
\dotsb \to X_{p-1}$.  These functors satisfy the following
formulas:
\[
\partial_{i}c_{p}=
\begin{cases}
\id&i=0\\
c_{p-1}\partial_{i-1}&i>0
\end{cases}
\qquad\qquad 
s_{i}c_{p}=
\begin{cases}
c_{p+1} c_{p}&i=0\\
c_{p+1}s_{i-1}&i>0
\end{cases}
\]
It follows that the functors $c\subdot$ specify a simplicial
contraction of $uM\subdot Z$, i.e., a retraction $C(uM\subdot Z)\to
uM\subdot Z$ of the inclusion $uM\subdot Z\to C(uM\subdot Z)$, where
$C(uM\subdot Z)$ is the join of a point with $uM\subdot Z$
(cf. \cite[5.1,6.2]{May}).

\section{The Proof of Theorems~\ref{thmpizero} and~\ref{thmemreduc}}\label{secpizero}

In this section we prove Theorem~\ref{thmpizero} from the introduction
and Theorem~\ref{thmemreduc} from the previous section, completing the
proof of Theorem~\ref{thmbigthm}.  The heart of the argument is
Lemma~\ref{lemmainlem} below that compares the nerves of topologized
categories with the nerves of the discrete categories obtained from
forgetting the topology.  This allows us to convert nerves constructed
out of sets of maps to nerves constructed out of the corresponding
spaces of maps.  This facilitates the comparison of the categories of
$R$-modules with the corresponding algebraic categories in
Theorem~\ref{thmpizero} and Theorem~\ref{thmemreduc} because, when the
domain is nice (homotopy equivalent to a cell $R$-module), the space
of maps between objects of $\aP_{0}^{0}$ (in the case of the
Theorem~\ref{thmpizero}) and the space of maps from a $K(\pi,n+1)$
$R$-module to an object in $\aP_{0}^{n+1}$ (in the case of
Theorem~\ref{thmemreduc}) are homotopy discrete with components the
appropriate corresponding sets of maps of $\pi_{0}R$-modules.

Because the hypotheses for Lemma~\ref{lemmainlem} are complicated, we
state it only for the categories $F'_{q}\aP$ (and $wF'_{q}\aP$), but
the proof uses few of the particulars of these categories and works
very generally for tensored topological categories.  For the
statement, let $B$ denote a fixed object of $F'_{q}\aP$, and
consider the over-category $F'_{q}\aP/B$: An object is a map $\alpha
\colon A\to B$ in $F'_{q}\aP$, and a morphism is a commutative
triangle.  A morphism $\alpha \to\alpha'$ between objects in
$F'_{q}\aP/B$ is a weak equivalence if the underlying map $A\to A'$ in
$F'_{q}\aP$ is a weak equivalence; let $w(F'_{q}\aP/B)$ denote the
subcategory of weak equivalences.  Let $\aD$ be a full subcategory of
$w(F'_{q}\aP/B)$ such that:
\begin{enumerate}
\item If $\alpha_{0} \colon A\to B$ is an object of $\aD$, and 
$\alpha_{1}\colon A\to B$ is homotopic to $\alpha_{0}$ in $F'_{q}\aP$, then
$\alpha_{1}$ is also an object of $\aD$
\item If $\alpha \colon A\to B$ is an object of $\aD$, then $A\sma I_{+}\to B$ is an
object of $\aD$.
\end{enumerate}
In (ii) we have in mind the map obtained by composing the
projection $A\sma I_{+}\to A$ with $\alpha$, but by (i), any map that restricts at
some point on $I$ to a map homotopic to $\alpha \colon A\to B$ is
an object of $\aD$.  Note that since $I^{n}$ is homeomorphic to
$\Delta[n]$, when $A\to B$ is an object of $\aD$ we can regard $A\sma
\Delta[n]_{+}\to B$ as an object of $\aD$.

We regard $\aD$ as a discrete category, but there is a natural
topology on both the objects and morphisms of $\aD$; let
$\aD^{\sing}\subdot$ denote the simplicial category obtained by taking
the total singular complex of the objects and morphisms.  More
concretely, an object of $\aD^{\sing}_{n}$ is a map $\alpha \colon
A\sma \Delta[n]_{+}\to B$ in $F'_{q}\aP$ that at each point of
$\Delta[n]$ restricts to an object of $\aD$ (or, equivalently, that is
itself an object of $\aD$).  A map $f\colon \alpha \to \alpha'$ in
$\aD^{\sing}_{n}$ is a map $f\colon A\sma \Delta[n]_{+}\to A'$ in
$wF'_{q}\aP$ such that the composite
\[
\xymatrix@C=4em{%
A\sma \Delta[n]_{+}\ar[r]^(.4){\id_{A} \sma \triangle}
&A\sma \Delta[n]_{+}\sma \Delta[n]_{+}\ar[r]^(.6){f\sma \id_{\Delta[n]}}
&A' \sma \Delta[n]_{+}\ar[r]^(.6){\alpha'}
&B
}
\]
is $\alpha$.  The category $\aD^{\sing}_{0}$ is $\aD$, and so we get a
simplicial functor $\aD\to \aD^{\sing}\subdot$, regarding $\aD$ as a
constant simplicial category.

\begin{lem}\label{lemmainlem}
With hypotheses and notation above, the inclusion $\aD\to
\aD^{\sing}\subdot$ induces a weak equivalence on nerves.
\end{lem}

\begin{proof}
It suffices to show that the iterated degeneracy functor
$s_{0}^{n}\colon \aD\to \aD^{\sing}_{n}$ induces a weak equivalence on
nerves.  This functor has a left adjoint $U$ that takes an object
$\alpha \colon A\sma
\Delta[n]_{+}\to B$ in $\aD^{\sing}_{n}$ and regards it as an object
of $\aD$ and takes a 
map $f\colon \alpha \to \alpha'$ in $\aD^{\sing}_{n}$ to the composite
\[
\xymatrix@C=4em{%
A\sma \Delta[n]_{+}\ar[r]^(.4){\id_{A} \sma \triangle}
&A\sma \Delta[n]_{+}\sma \Delta[n]_{+}\ar[r]^(.6){f\sma \id_{\Delta[n]}}
&A' \sma \Delta[n]_{+},
}
\]
a map in $\aD$.  The unit $\Id\to s_{0}^{n} \circ U$ in
$\aD^{\sing}_{n}$ is induced by the identity map in $\aD$: The identity on
$A\sma \Delta[n]_{+}$ in $\aD$ specifies a map from $\alpha\colon A\sma
\Delta[n]_{+}\to B$ to
\[
s_{0}^{n}U\alpha \colon (A\sma \Delta[n]_{+})\sma \Delta[n]_{+}\to B
\]
in $\aD^{\sing}_{n}$. The counit $U \circ s_{0}^{n} \to \Id$ in $\aD$ is
induced by the projection $(-)\sma \Delta[n]_{+}\to (-)$.
\end{proof}

We now move on to the proof of Theorem~\ref{thmpizero}.  Recall that
$\aE$ is the exact category of finitely-generated left
$\pi_{0}R$-modules.  Since the forgetful functor from
$w\Sdot[p]\aP_{0}^{0}$ to $wF_{p-1}\aP_{0}^{0}$ and the forgetful
functor from $w\Sdot[p]\aE$ to $wF_{p-1}\aE$ are equivalences of
categories, it suffices to show that the functor $\pi_{0}\colon
wF_{p-1}\aP_{0}^{0}\to wF_{p-1}\aE$ induces a weak equivalence on
nerves.  Let $\aC_{0}^{0}$ be the full subcategory of $\aP_{0}^{0}$ of
objects homotopy equivalent to cell $R$-modules; then $\aC_{0}^{0}$ is
a Waldhausen subcategory of $\aP_{0}^{0}$ and $F_{p-1}\aC_{0}^{0}$ is
a Waldhausen subcategory of $F_{p-1}\aP_{0}^{0}$.  The following lemma
allows us to work with $F_{p-1}\aC_{0}^{0}$ in place of
$F_{p-1}\aP_{0}^{0}$.

\begin{lem}\label{lemapplymainlem}
The inclusion $wF_{p-1}\aC_{0}^{0}\to wF_{p-1}\aP_{0}^{0}$ induces a
weak equivalence of nerves.
\end{lem}

\begin{proof}
Applying Quillen's Theorem~A \cite{QuillenAK}, it suffices to show
that for every object $B$ in $wF_{p-1}\aP_{0}^{0}$, the comma category
$\aD=wF_{p-1}\aC_{0}^{0} \downarrow B$ is weakly contractible.
Since the functor in question is the inclusion of a full subcategory,
the comma category $\aD$ is the full subcategory of
$w(F_{p-1}\aP_{0}^{0}/B)$ of maps 
$A\to B$ where $A$ is an object in $wF_{p-1}\aC_{0}^{0}$, and so is
therefore also a full subcategory of $w(F'_{p-1}\aP/B)$.
Applying Lemma~\ref{lemmainlem}, it suffices to show that the nerve of
$\aD^{\sing}\subdot$ is weakly contractible.  Consider the functor
$\pi_{0}$ from $\aD^{\sing}\subdot$ to the over-category
$wF_{p-1}\aE/\pi_{0}B$.  Because the constituent $R$-module $A_{i}$ of
an object $A$ of $F_{p-1}\aC_{0}^{0}$ is homotopy equivalent to a cell
$R$-module, the space of $R$-module maps $A_{i}\to B_{i}$ has the
correct homotopy type and so is homotopy discrete with components
$\Hom_{Z}(\pi_{0}A_{i},\pi_{0}B_{i})$.  Moreover, since the
constituent maps of $A$ are Hurewicz cofibrations, the space of maps
$A\to B$ in $wF_{p-1}\aP_{0}^{0}$ is likewise homotopy discrete and
with components $wF_{p-1}\aE(\pi_{0}A,\pi_{0}B)$; the natural topology
on the set of objects of $wF_{p-1}\aC_{0}^{0} \downarrow B$ is
homeomorphic to the disjoint union over the objects $A$ of
$F_{p-1}\aC_{0}^{0}$ of these spaces.  The natural topology on the set
of all morphisms in $wF_{p-1}\aC_{0}^{0} \downarrow B$ decomposes into
a disjoint union over pairs of objects $A,A'\in F_{p-1}\aC_{0}^{0}$ of
the spaces $wF_{p-1}\aP(A,A')\times wF_{p-1}\aP(A',B)$ with their
natural topology.  Again this is homotopy discrete with components the
disjoint union of $wF_{p-1}\aE(\pi_{0}A,\pi_{0}A')\times
wF_{p-1}\aE(\pi_{0}A',\pi_{0}B)$.  We conclude that the category
$\aD^{\sing}\subdot$ has a homotopy discrete simplicial set of
objects, has a homotopy discrete simplicial set of morphisms, has
category of components equivalent to $wF_{p-1}\aE/\pi_{0}B$, and the
induced map on nerves is a weak equivalence.  Since
$wF_{p-1}\aE/\pi_{0}B$ has a final object (namely, $\pi_{0}B$), it has
a contractible nerve.
\end{proof}

We can also apply Lemma~\ref{lemmainlem} with $\aD=wF_{p-1}\aC_{0}^{0}$
(with $B$ the final object).  The simplicial category
$\aD^{\sing}\subdot$ then has a homotopy discrete simplicial set of
objects, has a homotopy discrete simplicial set of maps, has category
of components equivalent to $wF_{p-1}\aE$, and has nerve weakly
equivalent to the nerve of $wF_{p-1}\aE$.  Since the composite functor
$wF_{p-1}\aC_{0}^{0}\to wF_{p-1}\aE$ is $\pi_{0}$, it follows that
$\pi_{0}\colon wF_{p-1}\aP_{0}^{0}\to wF_{p-1}\aE$ induces a weak
equivalence of nerves, and this completes the proof of
Theorem~\ref{thmpizero}.

The proof of Theorem~\ref{thmemreduc} is only slightly more
complicated.  By Proposition~\ref{propuff} in the previous section, it
suffices to show that the functor $\pi_{n+1}\colon
uF^{f}_{p-1}\aP_{0}^{n+1}\to uM_{p}Z$ induces a weak equivalence on
nerves.  Now let $C_{n+1}^{n+1}F^{f}_{p-1}$ be the full subcategory of
$uF^{f}_{p-1}\aP_{0}^{n+1}$ consisting of those sequences $A_{0}\to
\dotsb \to A_{p-1}$ where each map $A_{i}\to A_{i+1}$ is a Hurewicz
cofibration, each $A_{i}$ is homotopy equivalent to a cell $R$-module,
and $\pi_{q}A_{i}=0$ for $q\neq n+1$.  In other words, the objects of
$C_{n+1}^{n+1}F^{f}_{p-1}$ are the sequences of Hurewicz cofibrations
between homotopy cell $K(\pi,n+1)$ $R$-modules.  In particular, any
map between objects of $C_{n+1}^{n+1}F^{f}_{p-1}$ is a weak
equivalence, and so we can also consider $C_{n+1}^{n+1}F^{f}_{p-1}$ as
a full subcategory of $wF_{p-1}\aP$.  Then as above, we have the
following lemma.

\begin{lem}\label{lempost}
The inclusion $C_{n+1}^{n+1}F^{f}_{p-1}\to uF^{f}_{p-1}\aP_{0}^{n+1}$ induces a
weak equivalence of nerves.
\end{lem}

The proof is essentially the same combination of
Lemma~\ref{lemmainlem} with Quillen's Theorem~A as in the proof of
Lemma~\ref{lemapplymainlem} above, using $M_{p}Z$ in place of
$F_{p-1}\aE$: For an object $B$ in $uF^{f}_{p-1}\aP_{0}^{n+1}$, the
space of $R$-module maps from a homotopy cell $K(\pi,n+1)$ $R$-module
into the $R$-module $B_{i}$ (which has zero homotopy groups above
degree $n+1$) is homotopy discrete with components
$\Hom_{Z}(\pi,\pi_{n+1}B_{i})$ and for any object $A$ of
$C_{n+1}^{n+1}F^{f}_{p-1}$, the space of maps $A\to B$ in
$uF_{p-1}\aP_{0}^{n+1}$ is homotopy discrete with components
$uM_{p}Z(\pi_{n+1}A,\pi_{n+1}B)$.  As in the proof of
Lemma~\ref{lemapplymainlem}, it follows that the functor $\pi_{n+1}$
induces a weak equivalence of nerves from
$C_{n+1}^{n+1}F^{f}_{p-1} \downarrow B$ to $uM_{p}Z/\pi_{n+1}B$, which is
contractible.

Applying Lemma~\ref{lemmainlem} again with $\aD=C_{n+1}^{n+1}F^{f}_{p-1}$
and $B$ the final object, we see that the functor $\pi_{n+1}$ from
$C_{n+1}^{n+1}F^{f}_{p-1}$ to $uM_{p}Z$ induces a weak equivalence of
nerves.  It follows that the functor $\pi_{n+1}$ from
$uF^{f}_{p-1}\aP_{0}^{n+1}$ to $uM_{p}Z$ induces a weak equivalence of
nerves, and this completes the proof of Theorem~\ref{thmemreduc}.


\bibliographystyle{plain}

\end{document}